\documentclass[a4paper,11pt,english]{article}

\usepackage{amsmath,amsfonts,amssymb,amsthm}
\usepackage[all]{xy}
\usepackage{geometry,babel}

\newtheorem {thm}{Theorem}
\newtheorem* {thm*}{Theorem}

\newtheorem {cor}[thm]{Corollary}

\newtheorem* {cor*}{Corollary}
\newtheorem {lem}[thm]{Lemma}
\newtheorem {prop}[thm]{Proposition}

\newtheorem {rem}[thm]{Remark}

\theoremstyle{definition}

\newtheorem* {conj*}{Conjecture}

\newtheorem* {quest*}{Question}
\newtheorem* {problem}{Open problem}

\numberwithin{thm}{section}

\DeclareMathOperator{\End}{End}

\DeclareMathOperator{\Gal}{Gal}

\DeclareMathOperator{\Frob}{Frob}

\DeclareMathOperator{\SL}{SL}
\DeclareMathOperator{\GL}{GL}

\DeclareMathOperator{\Tr}{Tr}
\DeclareMathOperator{\Aut}{Aut}

\newcommand{\C}{\mathbb{C}}

\newcommand{\F}{\mathbb{F}}
\newcommand{\G}{\Gamma}
\newcommand{\N}{\mathbb{N}}
\newcommand{\Q}{\mathbb{Q}}
\newcommand{\Z}{\mathbb{Z}}

\renewcommand{\L}{\Lambda}

\newcommand{\g}{\gamma}
\newcommand{\e}{\varepsilon}
\newcommand{\p}{\mathfrak{p}}
\newcommand{\q}{\mathfrak{q}}

\newcommand{\z}{\zeta}

\renewcommand{\a}{\alpha}
\renewcommand{\b}{\beta}
\renewcommand{\l}{\lambda}

\newcommand{\Kbar}{\bar{K}}
\newcommand{\Fbar}{\bar{\F}}
\newcommand{\abar}{\bar{\a}}
\newcommand{\lbar}{{\bar\l}}
\newcommand{\zbar}{\bar\z}

\newcommand{\ee}{{\,\approx\,}}
\newcommand{\nee}{{\,\not\approx\,}}

\newcommand{\nl}{{v_\ell\,}}

\author{Chris Hall and Antonella Perucca}
\title{Radical Characterizations of Elliptic Curves}
\date{}

\begin{document}

\maketitle


\begin{abstract}
Let $K$ be a number field, and let $E$ be an elliptic curve over $K$. A famous result by Faltings of 1983 can be reformulated for elliptic curves as follows: if $S$ is a set of primes of good reduction for $E$ having density one, then the $K$-isogeny class of $E$ is determined by the function $\p\rightarrow \#E(k_{\p})$, which maps a prime in $S$ to the size of the group of points over the residue field.  In this paper, we prove that it suffices to look at the radical of the size.
\end{abstract}


\section{Introduction}

Let $A,A'$ be abelian varieties over a number field $K$, and let $S$ be the set of common finite primes $\p\subset K$ of good reduction.  A well-known theorem of Faltings implies that $A,A'$ are $K$-isogenous if and only if they have the same $L$-series (cf.~proposition~\ref{prop:faltings}).  The $L$-series of $A$ is determined, in part, by the function $\nu:\p\in S\mapsto\#A(k_\p)$, and in this paper we consider other functions which one can use to characterize $K$-isogeny.

If $A,A'$ are elliptic curves, then Faltings's theorem implies $A,A'$ have the same $L$-series if and only if $\# A(k_\p)=\#A'(k_\p)$ for every $\p\in S$ (cf.~lemma~\ref{nu-isogenous}).  A weaker condition that one could ask for is that these integers have the same radical, that is, $\ell\mid\#A(k_\p)$ if and only if $\ell\mid\# A'(k_\p)$, for every prime $\ell$ and every $\p\in S$.

\begin{thm}\label{thm1}
Suppose $S'\subseteq S$ has density one and $\L\subseteq\N$ is an infinite set of primes.  If $A,A'$ are elliptic curves over $K$, then the following are equivalent:
\begin{enumerate}
\item $A,A'$ are $K$-isogenous;
\item $\ell\mid\#A(k_\p)$ if and only if $\ell\mid\#A'(k_\p)$, for every $\ell\in\L$ and $\p\in S'$.
\end{enumerate}
\end{thm}

We prove theorem~\ref{thm1} in section~\ref{sec:thm1:proof}.  The proof makes crucial use of the theory developed by Serre in \cite{Serre} as well as the results of Frey and Jarden \cite{FreyJarden} on the Galois modules $A[\ell]$.  If $A$ is an elliptic curve and $A'=A\times A$, then 2.~holds but 1.~does not.  Nonetheless, the strategy behind our proof seems reasonable for pairs of square free abelian varieties, that is, abelian varieties $A,A'$ each of which whose isogeny factors are distinct.  We assume $A,A'$ are elliptic curves because our proof uses explicit knowledge of the structure of the image of the $\ell$-adic representation for $\ell\gg 0$.
A possible direction of further research is studying the following:

\begin{problem}
Let $A,A'$ be square-free abelian varieties defined over a number field, which are non-isogenous. Let $m,m'\geq 0$. Determine whether the following set has a positive density for $\ell\gg0$:
$$S_{\ell}=\{\p\in K: \nl(\#A(k_\p))=m,\, \nl(\#A'(k_\p))=m' \}\,.$$
\end{problem}

\subsection{Notation}

The following notation occurs frequently throughout the paper:

\begin{itemize}
\item $\ell$ : rational prime
\item $K$ : a number field
\item $G_K$ : absolute Galois group of $K$
\item $A,A'$ : abelian varieties (over $K$)
\item $E = \End(A)\otimes\Q$
\item $K_{\ell}$ : splitting field of the $\ell$-torsion of $A$
\item $L/K$ : a finite extension
\item $S_L$ : set of non-archimedean primes $\p\subseteq L$
\item $S_L(A)\subseteq S_L$ : primes of good reduction for $A$
\item $S_L(A,A')$ : common primes of good reduction for $A,A'$
\item $G_S$ : Galois group of the maximal extension in $\bar{K}$ unramified over $S\subseteq S_K$
\end{itemize}
Unless explicitly stated otherwise, we assume all abelian varieties and morphisms are defined over $K$.


\section{Arithmetic of elliptic curves}\label{sec:var_fibers}


\subsection{Local $L$-functions}\label{sec:locall}

Suppose $A$ is an elliptic curve and $\p\in S_K(A)$, and consider the following definitions
\[
	q_\p := \#k_\p,
	\quad
	a_\p := 1-\#A(k_\p)+q_\p,
	\quad
	\L_\p(T,A) := 1 - a_\p T + q_\p T^2.
\]
The special fiber $A_\p$ is an elliptic curve over $k_\p$, and $\L_\p(T,A)$ is the numerator of its so-called Hasse-Weil zeta function.

We call $\L_\p(T,A)$ the {\it local $L$-function} of $A$ at $\p$.  One basic property it has is that it satisfies the following identity:
\[
	\L_\p(1,A) = 1 - a_\p + q_\p = \# A(k_\p).
\]
A much deeper property $\L_\p(T,A)$ satisfies, the so-called Riemann hypothesis, is that the reciprocals $\a_\p,\b_\p\in\C^\times$ of its zeros satisfy $|\a_\p|=|\b_\p|=\sqrt{q_\p}$ (cf.~\cite[ch.~V, th.~2.3.1]{Silverman1new}).  By using the quadratic formula one can verify the following equivalences:
\[
	|\a_\p|=|\b_\p|=\sqrt{q_\p}
	\quad\Longleftrightarrow\quad
	\abar_\p = \b_\p
	\quad\Longleftrightarrow\quad
	|a_\p|\leq 2\sqrt{q_\p}.
\]
We suppose without loss of generality that $\a_\p$ has non-negative imaginary part so that it is well defined.

Recall that the absolute endomorphism ring of $A_\p$ is an order in the algebra $E_\p=\End(A_\p)\otimes\Q$.   Moreover, either $E_\p$ is a quadratic imaginary field or it is a quaternion algebra (cf.~\cite[ch.~V, th.~3.1]{Silverman1new}), and one says $A_\p$ is {\it ordinary} or {\it supersingular} respectively.

\begin{prop}\label{prop:waterhouse}
Suppose $\p\in S_K(A)$  does not lie over 2 or 3.  Then either $A_\p$ is ordinary and $[E_\p:\Q]=2$ or $A_\p$ is supersingular and exactly one of the following holds:
\begin{enumerate}
\item $[E_\p:\Q]=4$ and $\L_\p(T)\in\{\,(1\pm\sqrt{q_\p}\,T)^2\,\}$;
\item $[E_\p:\Q]=2$ and $\L_\p(T)\in\{\,1\pm\sqrt{q_\p}\,T+q_\p T^2\,\}$;
\item $[E_\p:\Q]=2$ and $\L_\p(T)=1+q_\p T^2$.
\end{enumerate}
Moreover, $A_\p$ is ordinary if and only if $\p\nmid a_\p$.
\end{prop}

\begin{proof}
This follows from \cite[th.~4.1]{Waterhouse}.
\end{proof}

Recall $E=\End(A)\otimes\Q$.

\begin{cor}\label{cor:Esplits}
Suppose $[E:\Q]=2$ and $E\subseteq K$.  If $\p\in S_K(A)$ does not lie over 2 or 3, then $E$ splits $\L_\p(T)$. 
\end{cor}

\begin{proof}
There is an embedding $E\to E_\p$, and thus proposition~\ref{prop:waterhouse} implies either $E=E_\p$ or $[E_\p:\Q]=4$ and $\L_\p(T)$ splits over $\Q$.  The Cayley-Hamilton theorem implies $\L_\p(\Frob_\p,A)=0$ in $E_\p$, and thus if $E=E_\p$, then $\L_p(T,A)$ has a zero in $E$ and so splits in $E$.
\end{proof}


\subsection{Cartan Subgroups}

Suppose $\ell$ is odd, and let $C_\ell\subseteq\GL_2(\F_\ell)$ be a \emph{Cartan subgroup}.  This means that, up to conjugating by some $g\in\GL_2(\F_\ell)$, for some $d\in\F_\ell^\times$ we have 
\[
	C_\ell = \left\{
		\left(
		\begin{array}{cc}
		a & b \\
		bd & a
		\end{array}
		\right)
		\ : \ 
		a,b\in\F_\ell
	\ \right\}\cap\GL_2(\F_\ell).
\]
The elements with $b=0$ correspond to the scalars $\F_\ell^\times\subseteq\GL_2(\F_\ell)$, and if $d=\delta^2$ with $\delta\in\Fbar_\ell$, then the $\Fbar_\ell$-eigenspaces of the non-scalar elements have the form $y=\pm\delta x$, with eigenvalues $a\pm b\delta$.

If $\delta\in\F_\ell$, then $C_\ell$ is called a \emph{split} Cartan subgroup and it consists of all elements which are diagonal with respect to a fixed basis of $\F_\ell^2$. Thus there is an isomorphism $C_\ell\simeq\F_\ell^\times\times\F_\ell^\times$, defined up to post-composition with $(x,y)\mapsto (y,x)$, and we have $\#C_\ell=(\ell-1)^2$.

If $\delta\not\in\F_\ell$, then $[\F_\ell(\delta):\F_\ell]=2$ and $C_\ell$ is called a \emph{non-split} Cartan subgroup.  The elements of $C_\ell$, plus the zero matrix, form a commutative division subring of $M_2(\F_\ell)$ with $\ell^2$ element.  Thus there is an isomorphism $C_\ell\simeq\F_{\ell^2}^\times$, defined up to post-composition with $x\mapsto x^\ell$, and we have $\#C_\ell=\ell^2-1$.


\subsection{Galois Theory}\label{sec:galois}

Let $K_\ell=K(A[\ell])$, and let $G_\ell$ be the Galois group of $K_\ell/K$.  We choose a basis of $A[\ell]$ so that we obtain an isomorphism $\Aut(A[\ell])\simeq\GL_{2g}(\F_\ell)$, for $g=\dim(A)$, and thus there is a natural embedding $G_\ell\to\GL_{2g}(\F_\ell)$ defined up to conjugacy.

We fix a polarization of $A$ and suppose $\ell$ does not divide its degree so that one can define the Weil pairing on $A[\ell]$.  The pairing takes its values in $\mu_\ell$, the group of $\ell$th roots of unity, so its existence implies $\mu_\ell\subseteq K_\ell$.  It is also Galois equivariant, that is, the following identity holds for all $\g\in G_\ell$:
\begin{equation}\label{eqn:weil}
	\langle P^\g,Q^\g\rangle = \langle P,Q\rangle^\g,
	\quad
	\forall P,Q\in A[\ell].
\end{equation}

We write $H_\ell\subseteq G_\ell$ for the Galois group of $K_\ell/K(\mu_\ell)$.  There is a natural embedding $G_\ell/H_\ell\to\Aut(\mu_\ell)=\F_\ell^\times$, and we write $\chi_\ell:G_\ell\to\F_\ell^\times$ for the composition of this embedding with the quotient map $G_\ell\to G_\ell/H_\ell$.

\begin{rem}
If $g=1$, then $\chi_\ell$ is the restriction to $G_\ell$ of $\det:\GL_2(\F_\ell)\to\F_\ell^\times$.
\end{rem}

\begin{rem}\label{rem:cyclotomic}
The induced homomorphism $\chi_\ell:G_K\to\F_\ell^\times$ is the cyclotomic character.
\end{rem}

Let $E=\End(A)\otimes\Q$.

\begin{prop}\label{prop:serre}
Suppose $A$ is an elliptic curve.  If $E\subseteq K$, then one of the following holds for all $\ell\gg 0$:
\begin{enumerate}
\item\label{case:Z} $E=\Q$ and $G_\ell=\GL_2(\F_\ell)$;
\item\label{case:cmR} $E\neq\Q$ and $G_\ell=C_\ell$ for some Cartan subgroup $C_\ell\subseteq\GL_2(\F_\ell)$.
\end{enumerate}
\end{prop}

\begin{proof}
The first part follows from \cite[th\'eor\`eme 2]{Serre}.  The second follows from \cite[\S 4.5, corollaire]{Serre} (cf.~section~\ref{sec:cm}).
\end{proof}

The following is a useful independence result:

\begin{prop}\label{prop:serre:indep}
Suppose $A$ is an elliptic curve and $L/K$ is a finite extension, and let $K'=KE$.  If $\ell\gg 0$, then $L\cap K_\ell\subseteq K'$.
\end{prop}

\begin{proof}
If we write $K'_\ell=K'(A[\ell])$, then $L\cap K'_\ell=K'$ implies $L\cap K_\ell\subseteq K'$.  Thus it suffices to show the $L\cap K'_\ell=K'$ for all $\ell\gg 0$, and so we replace $K,L$ by $K',LE$ and suppose $E\subseteq K$.  There are only finitely many possibilities for the intersection $L\cap K_\ell$ since $L$ is a finite (separable) extension of $K$, so some extension $K''/K$ occurs for infinitely many $\ell$.  In particular, \cite[th.~3 and \S 4.5, cor.]{Serre} implies $K_{\ell_1}\cap K_{\ell_2}=K$ for $\ell_1>\ell_2\gg 0$, so the only possibility is $K''=K$.  That is, $L\cap K_\ell=K$ for every $\ell\gg 0$.
\end{proof}

Suppose $A'$ is another polarized abelian variety, and let $K'_\ell$, $G'_\ell$, and $\chi'_\ell$ be defined accordingly for $A'$.  Consider the compositum $K_\ell K_\ell'$ as an extension of $K$ and let $\G_\ell\subseteq G_\ell\times G'_\ell$ be its Galois group.

\begin{lem}\label{lem:compositum}
$\G_\ell$ lies in the subgroup of $(\g,\g')\in G_\ell\times G'_\ell$ satisfying $\chi_\ell(\g)=\chi'_\ell(\g')$.
\end{lem}

\begin{proof}
As noted in remark~\ref{rem:cyclotomic}, the induced characters $\chi_\ell:G_K\to\F_\ell^\times$ and $\chi'_\ell:G_K\to\F_\ell^\times$ are both the cyclotomic character, that is, they are the same character.  Therefore the compositions of the projections $\G_\ell\to G_\ell,\G_\ell\to G'_\ell$ with the respective restrictions of $\chi_\ell,\chi'_\ell$ are the same map $\G_\ell\to\F_\ell^\times$, that is, $\chi_\ell(\g)=\chi'_\ell(\g')$ for each $(\g,\g')\in\G_\ell$.
\end{proof}


\section{Complex Multiplication}\label{sec:cm}

Throughout this section we assume $A$ is an elliptic curve and that both $[E:\Q]=2$ and $E\subseteq K$.  Thus $E$ is a quadratic imaginary field and $\End_K(A)\otimes\Q=E$ (cf.~\cite[\S II.2, th.~2.2]{Silverman}).  Unless stated otherwise, we also suppose that $A,A'$ are $\bar{K}$-isogenous so that $\End(A')\otimes\Q=E$ and define $S=S_K(A,A')$.

Suppose $\ell\geq 3$ is unramified in $E$, and fix a prime $\l\in S_E$ lying over $\ell$ and let $\lbar\in S_E$ be its $\Gal(E/\Q)$-conjugate.  The modules $A[\l]$ and $A[\lbar]$ are well-defined, and we have the following isomorphism of Galois modules:
\begin{equation}\label{eqn:cm}
	A[\ell] \simeq
		\begin{cases}
		A[\l]\times A[\lbar] & \mbox{$\ell$ split in $E$} \\
		\ \quad A[\l] & \mbox{$\ell$ inert in $E$}
		\end{cases}.
\end{equation}

We note that $A[\l]$ is a $\F_\l[G_K]$-module satisfying $\dim_{\F_\l}(A[\l])=1$, thus $\Aut_{\F_\l}(A[\l])=\F_\l^\times$.  Moreover, if $\ell\gg 0$, then the projections $A[\ell]\to A[\l],A[\lbar]$ lift to endomorphisms of $A$.  In particular, for almost all $\ell$, the action of $G_\ell$ must commute with the projections $A[\ell]\to A[\l],A[\lbar]$ and respective actions of $\F_\l^\times,\F_\lbar^\times$.  The subgroup $C_\ell\subseteq\Aut(A[\ell])$ of all elements with this property is a Cartan subgroup.  More precisely, $C_\ell$ is split if and only if $\ell$ splits in $K$, and then $C_\ell=\Aut(A[\l])\times\Aut(A[\lbar])=\F_\l^\times\times\F_\lbar^\times$.  Otherwise $C_\ell$ is non-split and $\ell$ is inert in $K$, and then $C_\ell=\F_\l^\times$.

Let $K_\l=K(A[\l])$ and $G_\l\subseteq\Aut(A[\l])$ be the Galois group of $K_\l/K$.

\begin{prop}\label{cor:cmserre}
If $\l\gg 0$, then there is a canonical isomorphism $G_\l=\F_\l^\times$.
\end{prop}

\begin{proof}
Suppose $\ell$ is large so that proposition~\ref{prop:serre} implies $G_\ell=C_\ell$.  If $C_\ell$ is non-split, then $G_\l=G_\ell=C_\ell=\F_\l^\times$, so suppose $C_\ell$ is split and thus $C_\ell=\F_\l^\times\times\F_\lbar^\times$ and $K_\ell=K_\l K_{\lbar}$.  There are natural embeddings $G_\ell\subseteq G_\l\times G_{\lbar}$ and $G_\l\times G_{\lbar}\subseteq C_\ell$ and their composition is a bijection, hence  $G_\l\times G_\lbar=C_\ell$ and thus $G_\l=\F_\l^\times$.
\end{proof}

If we write $S_{\ell}\subseteq S$ for the subset of $\p$ not lying over $\ell$, then $K_\ell/K$ is unramified over $S_\ell$ and hence so is $K_\l/K$.  That is, the composed character $G_K\to G_\l\subseteq\F_\l^\times$ factors through $G_K\to G_{S_\ell}$, and we write $\psi_\l:G_{S_\ell}\to\F_\l^\times$ for the corresponding character.  For split $C_\ell$, let $\psi_\lbar:G_{S_\ell}\to\F_\lbar^\times$ be defined similarly, and otherwise let $\psi_\lbar:G_{S_\ell}\to\F_\l^\times$ be the composition of $\psi_\l$ and the $\ell$th power map $\F_\l^\times\to\F_\l^\times$.

\begin{lem}\label{lem:aq}
If $\p\in S_\ell$ and if $\phi_\p\in G_{S_\ell}$ is an element in the conjugacy class of the Frobenius, then $(\psi_\l+\psi_\lbar)(\phi_\p)\equiv a_\p\pmod\ell$ and $\psi_\l\psi_\lbar(\phi_\p)\equiv q_\p\pmod\ell$.
\end{lem}

\begin{proof}
As noted in remark~\ref{rem:cyclotomic}, the composition of $G_K\to G_{S_\ell}\to G_\ell$ and $\det:G_\ell\to\F_\ell^\times$ is the cyclotomic character $\chi_\ell$.  Moreover, $\phi_\p$ acts in $\bar{k}_\p$ as $x\mapsto x^{q_\p}$, thus $\chi_\ell(\phi_\p)\equiv q_\p\pmod\ell$.  Similarly, the restriction $\Tr:G_\ell\to\F_\ell$ is well defined and  $\Tr(\phi_\p)\equiv a_\p\pmod\ell$.

If $C_\ell$ is split, then $C_\ell=\F_\l^\times\times\F_\lbar^\times$ and $G_\ell=(\psi_\l(G_K),\psi_\lbar(G_K))\subseteq C_\ell$, and $\det,\Tr$ are respectively given by the maps $(x,y)\mapsto xy$ and $(x,y)\mapsto x+y$.  Otherwise, $C_\ell=\F_\l^\times$ and $G_\ell=\psi_\l(G_K)\subseteq C_\ell$, and $\det,\Tr$ are respectively the norm and trace maps $\F_\l^\times\to\F_\ell^\times$ and $\F_\l\to\F_\ell$.
\end{proof}

For each $\p\in S_K(A)$, corollary \ref{cor:Esplits} implies the reciprocals $\a_\p,\abar_\p$ of the zeros of $\L_\p(T)$ lie in $E$ (and thus in $K$), so for $\p$ not dividing $\ell$, lemma~\ref{lem:aq} implies we have the following identities:
\begin{equation}\label{eqn:modl}
	\L_\p(T)\equiv (1-\a_\p T)(1-\abar_\p T)
	\equiv (1-\psi_\l(\phi_\p)T)(1-\psi_\lbar(\phi_\p)T)\pmod\l.
\end{equation}
Suppose $A'$ is $\bar{K}$-isogenous to $A$.  Then a similar congruence holds for $\L'_\p(T)$ and $\p\in S_K(A')$, and in particular, we have the following identities for $\p\in S$:
\begin{equation}\label{eqn:linL}
	\{\a_\p,\abar_\p\} \equiv \{\psi_\l(\phi_\p),\psi_\lbar(\phi_\p)\},
	\quad
	\{\a'_\p,\abar'_\p\} \equiv \{\psi'_\l(\phi_\p),\psi'_\lbar(\phi_\p)\}
	\pmod\l.
\end{equation}

Let $\e_\l:G_{S_\ell}\to\F_\l^\times$ be the character $\e_\l:=\psi'_\l/\psi_\l$ so that $A[\l]\otimes\e_\l\simeq A'[\l]$.

\begin{lem}\label{lem:boundedorder}
If $A,A'$ are $\bar{K}$-isogenous, then for some $m\geq 1$ and all $\l\gg 0$, the order of $\e_\l$ is at most $m$.
\end{lem}

\begin{proof}
Suppose $L/K$ is a finite Galois extension and over which $A,A'$ are isogenous, and let $S_\ell'\subseteq S_\ell$ be the subset of $\p$ which are unramified in $L$ so that $\psi_\l,\psi'_\l,\e_\l$ all factor through characters $G_{S'_\ell}\to\F_\l^\times$.  For each $\p\in S'_\ell$, choose $\q\in S_L$ dividing $\p$ and let $\phi_\p,\phi_\q\in G_{S'_\ell}$ be the corresponding Frobenius elements.  Then we have the following identities:
$$
	\psi_\l(\phi_\p)^{[k_\q:k_\p]} \equiv \psi_\l(\phi_\q),
	\quad
	\psi'_\l(\phi_\p)^{[k_\q:k_\p]} \equiv \psi'_\l(\phi_\q)
	\pmod\l.
$$

Suppose $\phi:A\to A'$ is an $L$-isogeny.  If $\ell$ does not divide $\deg(\phi)$, then $\phi$ induces an isomorphism $A[\l]\simeq A'[\l]$ of $G_L$-modules and thus the restrictions $\psi_\l:G_L\to\F_\l^\times$ and $\psi'_\l:G_L\to\F_\l^\times$ are isomorphic.  Therefore, since $L/K$ is Galois and so $[k_\q:k_\p]\mid [L:K]$, we have the following identities:
\[
	\e_\l(\phi_\p)^{[L:K]}
	\equiv (\psi_\l(\phi_\p)/\psi'_\l(\phi_\p))^{[L:K]}
	\equiv 1
	\!\!\pmod\l.
\]
That is, $\e_\l$ has order dividing $[L:K]$ for every $\l\gg 0$.
\end{proof}

\begin{cor}\label{cor:G_S}
Suppose $A,A'$ are $\bar{K}$-isogenous, and let $S=S_K(A,A')$.  If $\l\gg 0$, then $\e_\l$ factors through the quotient $G_{S_\ell}\to G_S$.
\end{cor}

\begin{proof}
Let $m$ be as in lemma~\ref{lem:boundedorder}, and suppose $\l\gg 0$ so that the order of $\e_\l$ is at most $m$.  Since $\e_\l$ takes values in $\F_\l^\times$, it has order prime to $\ell$, and thus $\e_\l$ is at worst tamely ramified over $\p$ dividing $\ell$.  Therefore it suffices to show the restriction of $\e_\l$ to the tame inertia group $I_\ell\subseteq G_{S_\ell}$, defined up to conjugacy, is trivial (cf. \cite[\S 1.6]{Serre}).  Equivalently, it suffices to show that the corresponding invariant, which is an element in $\Q/\Z$ with denominator coprime to $\ell$, is zero (cf. \cite[\S 1.7]{Serre}).

If $\l\gg 0$, then $\ell$ is unramified in $K$ and thus the corollaries of \cite[prop.~11, prop.~12]{Serre} imply that the respective invariants $x,x'$ of $\psi_\l,\psi'_\l$ are elements of
$$
	X = \{\ 0,\ 1/(\ell-1),\ 1/(\ell^2-1),\ \ell/(\ell^2-1)\ \}.
$$
Swapping $A,A'$ if necessary, we suppose $x\geq x'$.  Then the invariant of $\e_\l=\psi_\l/\psi'_\l$ is $x-x'$ and thus is an element of
$$
	X' = \{\ 0, 1/(\ell-1),\ 1/(\ell^2-1),\ \ell/(\ell^2-1),\ \ell/(\ell-1),\ 1/(\ell+1) \ \}
$$

If $m_\l$ is the order of $\e_\l$, then $x-x'=a/b$ is a multiple of $1/m_\l$, that is, $b\mid m_\l$ when $\gcd(a,b)=1$.  Thus if $a\neq 0$, then $b$ is one of $\ell-1$, $\ell+1$, or $\ell^2-1$.  In particular, since $b\leq m_\l\leq m$, if $\ell>m+1$, then we must have $a=0$, that is, the invariant of $\e_\l$ is trivial for $\l\gg 0$.
\end{proof}

Let $\mu\subset E$ by the subgroup of roots of unity.  We define $\mu\to\F_\l^\times$ to be the homomorphism induced by reduction modulo $\l$; it is injective, if $\ell\geq 5$.  If $\l$ is split, we define $\mu\to\F_\lbar^\times$ to be reduction modulo $\lbar$.

\begin{cor}\label{cor:mu_values}
Suppose $A,A'$ are $\bar{K}$-isogenous.  If $\l\gg 0$, then $\e_\l$ takes values in $\mu$ and $\e_\l=\e_{\l'}$ for all $\l'$ in some infinite $\L\subseteq S_E$.
\end{cor}

\begin{proof}
Suppose $\l\gg 0$ so that lemma~\ref{lem:boundedorder} implies $\e_\l$ has order at most $m$ and corollary~\ref{cor:G_S} implies $\e_\l$ factors through $G_{S_\ell}\to G_S$.  Hermite's theorem implies the set of characters of $G_S$ of order at most $m$ is finite, hence up to excluding finitely many $\l$ we may suppose $\e_\l=\e_{\l'}$ for all $\l'$ in an infinite $\L\subseteq S_E$.

Suppose $\p\in S_\ell$.  Up to swapping $\a_\p,\abar_\p$ and up to swapping $\a'_\p,\abar'_\p$, we apply (\ref{eqn:linL}) and assume the following identities holds for any $\l'\in\L$ not dividing $\p$:
\[
	\psi_{\l'}(\phi_\p)\equiv\a_\p,
	\quad
	\psi'_{\l'}(\phi_\p)\equiv \a'_\p,
	\quad
	\e_{\l'}(\phi_\p)\equiv \a_\p/\a'_\p
	\pmod{\l'}.
\]
Since this identity holds for infinitely many $\l'$ and since $\e_\l(\phi_\p)=\z_\p=\e_{\l'}(\phi_\p)$ for some root of unity $\z_\p\in\bar{E}$, we must have $\z_\p=\a_\p/\a'_\p$ and thus $\z_\p\in\mu\subset E$.
\end{proof}

\begin{rem}\label{rem:cmtwist}
Suppose $\l\gg 0$, and let $\L\subseteq S_E$ be as in corollary~\ref{cor:mu_values} and $\chi:G_S\to\mu$ be such that $\chi=\e_\l$ for every $\l\in\L$.  Then $\psi_\l\chi=\psi'_\l$, or equivalently, $A[\l]\otimes\chi\simeq A'[\l]$ for every $\l\in\L$.
\end{rem}

Suppose $\ell\geq 5$, and consider the following embedding $\mu\subseteq C_\ell$: if $C_\ell$ is split, then it is the product embedding $\mu\to\F_\l^\times\times\F_\lbar^\times=C_\ell$ given by $\z\mapsto(\z,\z)$ ($=(\z,1/\z)\in\F_\ell^\times\times\F_\ell^\times$ via the field isomorphisms $\F_\l=\F_\ell=\F_\lbar$); otherwise it is the embedding $\mu\to\F_\l^\times=C_\ell$.  Thus for any character $\chi:G_K\to\mu$ we can consider the twist $A[\ell]\otimes\chi$ as an $\F_\ell[G_K]$-module.

\begin{lem}\label{lem:tofulll}
If $A[\l]\otimes\chi\simeq A'[\l]$, then $A[\ell]\otimes\chi\simeq A'[\ell]$.
\end{lem}

\begin{proof}
If $\ell$ is inert, then there is nothing to prove, so suppose $\ell$ is split and $A[\l]\otimes\chi\simeq A'[\l]$.  On one hand, $C_\ell=\F_\l^\times\times\F_\lbar^\times$ and the composition of $\chi:G_K\to\mu\subseteq C_\ell$ with the projection $C_\ell\to\F_\lbar^\times$ and bijection $\F_\lbar=\F_\l$ is $\chi^{-1}$, thus $A[\ell]\otimes\chi$ is isomorphic to the direct sum of $A[\l]\otimes\chi$ and $A[\lbar]\otimes(1/\chi)$.  On the other hand, since $\psi_\l\chi=\psi'_\l$, lemma~\ref{lem:aq} implies $\psi_\l\psi_{\bar\l}=\psi'_\l\psi'_{\bar\l}$, thus $\psi'_{\bar\l}=\psi_{\bar\l}/\chi$ and so $A[\lbar]\otimes(1/\chi)\simeq A'[\lbar]$.  Together these imply $A[\ell]\otimes\chi\simeq A'[\ell]$.
\end{proof}

\begin{cor}\label{cor:cmtwist}
If $A,A'$ are $\bar{K}$-isogenous, then there is an infinite $\L_\Q\subseteq S_\Q$ and a character $\chi:G_K\to\mu$ such that $A[\ell]\otimes\chi\simeq A'[\ell]$ for every $\ell\in\L_\Q$.
\end{cor}

\begin{proof}
Let $\L_\Q=\{\l\cap\Q:\l\in\L\}$ and combine remark~\ref{rem:cmtwist} and lemma~\ref{lem:tofulll}.
\end{proof}


\section{Fibers of elliptic curves}

\subsection{Essentially Equal Functions}

We say that a pair of functions $f,f'$ defined on a subset of $S_K$ are {\it essentially equal} if and only if they are equal on some density-one subset, and then we write $f\ee f'$.

Given abelian varieties $A,A'$, let $\l,\l'$ be the functions on $S_K(A)$ which maps $\p\in S_K(A)$ to the respective local $L$-function $\L_\p(T,A),\L_\p(T,A')\in\Z[T]$ (cf.~section~\ref{sec:locall}).

\begin{prop}\label{prop:faltings}
$\l\ee\l'$ if and only if $A,A'$ are $K$-isogenous.
\end{prop}
\begin{proof}
See \cite[Corollary 2]{Faltings83}.
\end{proof}

Given an abelian variety $A$ and finite extension $L/K$, we consider the function $\phi_L$ which maps $\p\in S_L(A)$ to the isomorphism class of the group of points $A(k_\p)$ on the special fiber.  Given a function $\psi$ on the set of groups modulo isomorphism, we write $\psi\phi_L$ for the composed function which maps $\p\in S_L(A)$ to $\psi(A(k_\p))$.  We define $\phi'_L$ and $\psi\phi'_L$ similarly given an additional abelian variety $A'$.

The following lemma demonstrates $\ee$ behaves well with respect to base change:

\begin{lem}\label{lem:basechange}
If $A,A'$ are abelian varieties and if $\psi\phi_K\ee \psi\phi_K'$, then $\psi\phi_L\ee\psi\phi_L'$.
\end{lem}
\begin{proof}
Let $S_K'\subseteq S_K$ and $S_L'\subseteq S_L$ be the respective subsets of primes of degree one.  They have Dirichlet density one.  If $\q\in S_L'$ and if $\p:=\q\cap K$, then the embedding $k_\p\to k_{\q}$ is surjective and
$$
	\psi\phi_L(\q) = \psi\phi_K(\p) 
		= \psi\phi_K'(\p) = \psi\phi_L'(\q),
$$
and thus $\psi\phi_L(\q)= \psi\phi_L'(\q)$ for all $\q$ in the density-one set $S_L'$.
\end{proof}

\begin{rem}
The converse does not hold: if $A,A'$ are not $K$-isogenous but are $L$-isomorphic, then $\phi_{A,K}\nee\phi_{A',K}$ and $\phi_{A,L}\ee\phi_{A',L}$. 
\end{rem}

For each $A$ and $\ell$, we write $\nu_{A,K}$ for the function $\p\in S_K(A)\mapsto\#A(k_\p)$ and regard it as the composition of $\phi_{A,L}$ and the counting function $G\mapsto\#G$.  A priori, one could have $\l_{A,K}\nee\l_{A',K}$ and yet still have $\nu_{A,K}\ee\nu_{A',K}$, but the lemma shows this does not occur for $\dim(A)=1$:

\begin{lem}\label{nu-isogenous}
If $A,A'$ are elliptic curves over $K$, then $\nu_{A,K}\ee\nu_{A',K}$ if and only if $A,A'$ are $K$-isogenous.
\end{lem}
\begin{proof}
If we write $\l=\l_{A,K}$ and define $\l',\nu,\nu'$ similarly, then we the following equivalences between identities imply $\l\ee\l'$ if and only if $\nu\ee\nu'$:
\[
	\l(\p) = \l'(\p)
	\quad\Longleftrightarrow\quad
	a_\p = a'_\p
	\quad\Longleftrightarrow\quad
	\nu(\p) = 1-a_\p+q_\p = 1-a'_\p+q_\p = \nu'(\p).
\]
The lemma follows if we apply proposition~\ref{prop:faltings} to deduce that $A,A'$ are $K$-isogenous if and only if $\l\ee\l'$.

\end{proof}


\subsection{Radicals}\label{sec:radicals}

Suppose $\ell$ is a prime, and let $g=\dim(A)$.  Recall that the radical of a finite group $G$ is the square free product of the primes dividing $\#G$.  We consider only the $\ell$-part of the radical and define $\rho_\ell$ to be the following map:
$$
	S_K(A)\to\{0,1\} : \p \mapsto \min\{1,v_\ell(\# A(k_\p))\};
$$
Recall that in section~\ref{sec:galois} we defined $K_\ell=K(A[\ell])$ and $G_\ell=\Gal(K_\ell/K)$.  The following lemma gives a Galois-theoretic way to analyze $\rho_\ell$:

\begin{lem}\label{lem:frob}
Suppose $\p\in S_K(A)$ does not ramify in $K_\ell$ and $\q\in S_{K_\ell}$ lies over $\p$.  If $\phi_\q\in G_\ell$ is the Frobenius of $\q$, then $\rho_\ell(\p)=1$ if and only if $\det(\phi_\q-1)=0$.
\end{lem}

\begin{proof}
The embedding $A(k_\p)\to A(k_{\q})$ identifies $A(k_\p)[\ell]$ with $\ker(\phi_q-1)\subseteq A[\ell]$, hence $\ell\mid\#A(k_\p)$ if and only if 1 is an eigenvalue of $\phi_\q$.
\end{proof}

Recall that we defined $K'_\ell$ and $G'_\ell$ for $A'$ and that $\G_\ell\subseteq G_\ell\times G'_\ell$ denotes the Galois group of the compositum $K_\ell K'_\ell/K$.  Let $\rho_\ell'$ be defined accordingly for $A'$.

\begin{lem}\label{cjh:lem1}
If $\rho_\ell\ee\rho_\ell'$, then $\det(\g-1),\det(\g'-1)$ are both zero or both non-zero for every $(\g,\g')\in\G_\ell$.
\end{lem}

\begin{proof}
Suppose $S\subseteq S_K(A,A')$ has Dirichlet density one and that $\rho_\ell|_S=\rho_\ell'|_S$.  Let $S'\subseteq S$ be the subset consisting of the primes $\p$ which are unramified in $K_\ell K'_\ell$ and whose Frobenius conjugacy class in $\G_\ell$ contains $(\g,\g')$.  The density of $S'$ is positive, and for each $\p\in S'$, lemma~\ref{lem:frob} implies the values $\rho_\ell(\p)$, $\rho_\ell'(\p)$ respectively identify whether or not $\det(\g-1),\det(\g'-1)$ are non-zero, and thus the hypothesis $\rho_\ell(\p)=\rho_\ell'(\p)$ implies the determinants are both zero or both non-zero.
\end{proof}


\section{Proof of Theorem~\ref{thm1}}\label{sec:thm1:proof}

Let $A,A'$ be elliptic curves over $K$ and $\rho_\ell,\rho'_\ell$ the maps defined in section~\ref{sec:radicals}.  In this section we prove the following result:

\bigskip\noindent
{\bf Theorem~\ref{thm1}. \it
Suppose $S'\subseteq S=S_K(A,A')$ has density one and $\L\subseteq\N$ is an infinite set of primes.  If $A,A'$ are elliptic curves over $K$, then the following are equivalent:
\begin{enumerate}
\item $A,A'$ are $K$-isogenous;
\item $\rho_\ell(\p)=\rho'_\ell(\p)$ for every $\ell\in\L$ and $\p\in S'$.
\end{enumerate}
}

\medskip
\noindent
The implication $1\Rightarrow 2$ follows from lemma~\ref{nu-isogenous}, so we prove $2\Rightarrow 1$.  The structure of the proof is as follows:
\begin{enumerate}

\item We reduce to the case $E,E'\subseteq K$ for $E=\End(A)\otimes\Q$ and $E'=\End(A')\otimes\Q$.

\item We show $K(A[\ell])=K(A'[\ell])$ for all $\ell$ in an infinite $\L'\subseteq\N$ and deduce $A,A'$ are $\Kbar$-isogenous and $E=E'$.

\item We construct a character $\chi$ such that $A[\ell]\otimes\chi\simeq A'[\ell]$ for almost all $\ell\in\L'$.

\item We prove $A[\ell]\simeq A'[\ell]$ for all almost $\ell\in\L'$ and deduce that $A,A'$ are $K$-isogenous.
\end{enumerate}

\noindent
The remainder of this section is broken into four pieces, one for each of these steps.


\subsubsection*{Step 1}

Lemma~\ref{lem:basechange} implies it suffices to prove the theorem after replacing $K$ by a finite extension, and two applications of the following lemma, one with $A,A'$ swapped, imply that it suffices to prove theorem~\ref{thm1} over $KEE'=(KE)E'$:

\begin{lem}\label{lem:KE-isogenous}
If $A,A'$ are $KE$-isogenous, then $E'=E$ and $A,A'$  are $K$-isogenous.
\end{lem}

\begin{proof}
The identity $\End(A)\otimes\Q=\End(A')\otimes\Q$ holds for any pair of $\bar{K}$-isogenous abelian varieties, thus $E=E'$.  Moreover, if $E\subseteq K$ then $A,A'$ are $K$-isogenous, so suppose $E\not\subseteq K$.

Let $S\subseteq S_K(A,A')$ be the density-one subset of $\p$ which have degree one and which neither ramify in $KE$ or lie over 2 or 3.  We will show that $a_\p=a'_\p$ for every $\p\in S$, and then lemma~\ref{nu-isogenous} implies $A,A'$ are $K$-isogenous, so suppose $\p\in S$.

If $\q\in S_{KE}$ is a prime lying over $\p$, then $a_{\q}=a_{\q}'$ since $A,A'$ are $KE$-isogenous.  If $\p$ splits in $KE$, then we have $a_\p=a_{\q}$ and $a'_\q=a'_\p$ since $k_\q=k_\p$, thus $a_\p=a'_\p$.  Otherwise, $\p\in S$ is inert, thus \cite[ch. 10 \S 4 theorem 10]{Lang} implies $A,A'$ have supersingular reduction over $\p$.  Moreover, since $q_\p$ is prime and thus not a square, proposition~\ref{prop:waterhouse} implies $a_\p=a'_\p=0$.  Therefore $a_\p=a'_\p$ for every $\p\in S$ as claimed.
\end{proof}


\subsubsection*{Step 2}

We use the notation of section~\ref{sec:galois}, thus $K_\ell=K(A[\ell])$, $G_\ell=\Gal(K_{\ell}/K)\subseteq\GL_2(\F_\ell)$, and $H_\ell=G_\ell\cap\SL_2(\F_\ell)$ is the stabilizer of $K(\zeta_\ell)$.  Similarly, we have corresponding objects $K'_\ell,G'_\ell,H'_\ell$.  Finally, $\G_\ell\subseteq G_\ell\times G'_\ell$ is the Galois group of $K_\ell K'_\ell/K$.

The kernels of the projections $\G_\ell\to G_\ell$ and $\G_\ell\to G'_\ell$ project onto normal subgroups of $G'_\ell$ and $G_\ell$ respectively.  For example, lemma~\ref{lem:compositum} implies the intersection of $\G_\ell$ with the normal subgroup $\SL_2(\F_\ell)\times\{1\}\subseteq\GL_2(\F_\ell)\times\GL_2(\F_\ell)$ is the kernel of $\G_\ell\to G'_\ell$ and it projects isomorphically onto a normal subgroup of $G_\ell$ contained in $H_\ell$.  Moreover, this kernel is trivial if and only if $K_\ell\subseteq K_\ell'$, thus both kernels are trivial if and only if $K_\ell=K_\ell'$.

\begin{lem}
Suppose $E,E'\subseteq K$.  If $K_\ell\neq K_\ell'$ and if $\ell\gg 0$, then $\rho_\ell\nee\rho'_\ell$.
\end{lem}

\begin{proof}
Suppose $K_\ell\neq K_\ell'$ and $\ell\gg 0$, and without loss of generality suppose the kernel of $\G_\ell\to G'_\ell$ is non-trivial and thus projects to a non-trivial normal subgroup of $G_\ell$ contained in $H_\ell$.  Since $E\subseteq K$ and $\ell$ is large, proposition~\ref{prop:serre} implies $G_\ell=\GL_2(\F_\ell)$ or $G_\ell=C_\ell$ for some Cartan subgroup $C_\ell\subseteq\GL_2(\F_\ell)$.

In the first case, the $g=-1$ lies in every non-trivial normal subgroup of $H_\ell=\SL_2(\F_\ell)$ (cf.~\cite[lem.~2.2]{FreyJarden}).  In the second case, every  $g\in H_\ell$ is semisimple and satisfies $\det(g)=1$, so either $g=1$ or $\det(g-1)\neq 0$.  Either way, we can find an element $(g,1)$ in the kernel satisfying $\det(g-1)\neq 0$, and thus lemma~\ref{cjh:lem1} implies $\rho_\ell\nee\rho'_\ell$.
\end{proof}

Since by assumption $\rho_\ell\ee\rho'_\ell$ for infinitely many $\ell$, we conclude that $K_\ell=K_\ell'$ for infinitely many $\ell$.

\begin{prop}\label{prop:fj}
If $K_\ell=K'_\ell$ for infinitely many $\ell$, then $A,A'$ are $\bar{K}$-isogenous and $E=E'$.
\end{prop}

\begin{proof}
The varieties $A,A'$ are $\bar{K}$-isogenous by \cite[theorem A]{FreyJarden}, and thus $E=E'$.
\end{proof}


\subsubsection*{Step 3}

Let $\mu\subset E^\times$ by the subgroup of roots of unity.  If $[E:\Q]=2$, then recall that, for $\ell\geq 5$, there is an embedding $\mu\to C_\ell\subseteq\Aut(A[\ell])$ (cf.~section~\ref{sec:cm}).

\begin{lem}\label{lem:twist}
Suppose $E\subseteq K$ and $K_\ell=K'_\ell$ for infinitely many $\ell$.  Then there exist an infinite $\L'\subseteq S_\Q$ and a character $\chi:G_K\to\mu$ such that $A[\ell]\otimes\chi\simeq A'[\ell]$ for all $\ell\in\L'$.
\end{lem}

\begin{proof}
If $E=\Q$ and if $\ell\gg 0$, then $G_\ell=G'_\ell\simeq\GL_2(\F_\ell)$ and \cite[lem.~8]{Serre} implies the result.  Otherwise we can apply proposition~\ref{prop:fj} to deduce that $A,A'$ are $\bar{K}$-isogenous and apply corollary~\ref{cor:cmtwist} to conclude.
\end{proof}

\begin{rem}
We may take $\L'\subseteq\L$ in lemma~\ref{lem:twist}.
\end{rem}


\subsubsection*{Step 4}

If $\chi:G_K\to\mu$  and $\L'\subseteq\L$ are as in lemma~\ref{lem:twist}, then the following lemma implies $\chi$ must be trivial and hence $A[\ell]\simeq A'[\ell]$ for almost all $\ell\in\L'$:

\begin{lem}\label{lem:twist2}
Suppose $E\subseteq K$ and $\chi:G_K\to\mu$ is non trivial.  If $A[\ell]\otimes\chi \simeq A'[\ell]$ and if $\ell\gg 0$, then $\rho_\ell\nee\rho'_\ell$.
\end{lem}

\begin{proof}
Let $L/K$ be the splitting field of $\chi$.  By construction, it is a Galois extension and $\chi$ identifies $\Gal(L/K)$ with the subgroup $\chi(G_K)\subseteq \mu$.

Let $S'_\ell\subseteq S$ be the subset of $\p$ which split completely in $K_\ell/K$.  If $\p\in S'_\ell$, then $A(k_\p)[\ell]=A[\ell]$ and also $\ell\mid \#k_\p-1$ since the existence of the Weil pairing implies $\mu_\ell\subseteq k_\p$.  Hence $\rho_\ell(\p)=1$ and we have the congruences
\[
	\L_\p(T) \equiv (1 - T)^2
	\pmod\ell.
\]

Suppose $\p\in S'$ and $\ell\geq 5$, and let $\z_\p=\chi(\phi_\p)$, $\zbar_\p=1/\z_\p$ and let $\a_\p,\abar_\p\in\Fbar_\ell$ be the reciprocals of the eigenvalues of the image of $\phi_\p$ in $G_\ell$.  Suppose $A[\ell]\otimes\chi\simeq A'[\ell]$.  If $\mu=\mu_2$, then $\zbar_\p=\z_\p$ and $\L'_\p(T)=\L_\p(\z_\p T)$, and otherwise, up to swapping $\a_\p,\abar_\p$ or $\a'_\p,\abar'_\p$, we have the following congruences for any $\l\in S_E$ dividing $\ell$:
\[
	\a'_\p \equiv \psi'_\l(\phi_\p) \equiv \psi_\l(\phi_\p)\chi(\phi_\p)
	\equiv \z_\p\a_\p,\quad \abar'_\p\equiv \zbar_\p\abar_\p \pmod{\l}.
\]
That is, for any $\mu$ and all $\l\in S_E$ dividing $\ell$, we have the following congruences (cf.~(\ref{eqn:modl})):
\[
	\L_\p(T)\equiv (1-\a_\p T)(1-\abar_\p T),
	\quad
	\L'_\p(T)\equiv (1-\z_\p\a_\p T)(1-\zbar_\p\abar_\p T)
	\pmod\l.
\]
Therefore, if moreover $\p\in S'_\ell$, then we have the following congruences:
\[
	\#A'(k_\p)=\L'_\p(1)\equiv (1-\z_\p)(1-\zbar_\p)\pmod\l.
\]
In particular, if $\z_\p\neq 1$, then last term is non zero and so $\rho'_\ell(\p)=0$.

To complete the proof we recall that $E\subseteq K$ and observe that proposition~\ref{prop:serre:indep} implies $L\cap K_\ell=K$ for $\ell\gg 0$, thus the subset $S''_\ell\subseteq S'_\ell$ of $\p$ such that $\z_\p\neq 1$ has positive density.
\end{proof}

In summary, $A[\ell]\simeq A'[\ell]$ for infinitely many $\ell$, hence \cite[prop.~1.4]{FreyJarden} implies $A,A'$ are $K$-isogenous.~~Q.E.D.


\vspace{0.6cm}

\noindent {\it Chris Hall,} University of Wyoming\\
\noindent \text{E-mail}: chall14@uwyo.edu\\

\noindent {\it Antonella Perucca,} Research Foundation - Flanders (FWO)

\noindent \text{E-mail}: antonellaperucca@gmail.com\\

\end{document}